\documentclass[10pt,reqno]{amsart}
\usepackage{amsmath}
\usepackage{amsfonts}
\usepackage{amssymb}
\usepackage{amsthm}
\usepackage{amscd}

\usepackage{a4wide}

\newcommand{\Section}[1]{\section{#1} \setcounter{equation}{0}}

\newcommand\Q{\mathbb{Q}}
\newcommand\R{\mathbb{R}}
\newcommand\Z{\mathbb{Z}}

\newcommand{\vv}[1]{{\mathbf{#1}}}

\newtheorem{thm}{Theorem}
\newtheorem{lem}{Lemma}[section]
\newtheorem{cor}{Corollary}

\theoremstyle{remark}

\theoremstyle{definition}

\def\hs{{\mathcal H}^s}

\begin{document}
\large
\title{Diophantine approximation on planar curves: the convergence theory}
\author{R.C. Vaughan}
\address{RCV: Department of Mathematics, McAllister Building,
Pennsylvania State University, University Park, PA 16802-6401,
U.S.A.} \email{rvaughan@math.psu.edu}
\author{S. Velani}
\address{SV: Department of Mathematics, University of York, Heslington, York, YO10
5DD, U.K.} \email{slv3@york.ac.uk }
\thanks{\\ 
 SV: Royal Society University Research Fellow.\\
 }

\begin{abstract}
The  convergence theory for the set of  simultaneously
$\psi$-approximable points lying on a planar curve is established.
Our results complement the divergence theory developed in
 \cite{1} and thereby completes the general metric
theory for planar curves.  
\\
\\
\\
\begin{center} {\large \em
Dedicated to Walter Hayman and Klaus Roth on their eightieth
 birthdays. } \end{center}
\end{abstract}


\maketitle


\Section{Introduction and Statement of Results}

\subsection{The motivation} \label{s:Intro}

In this paper we establish variants of Conjecture 1 of Beresnevich
{\it et al} \cite{1} that are sufficient to establish Conjecture 2
and Conjecture H of \cite{1}. Conjecture 1 is firmly rooted in
replacing  the upper bound in Huxley's theorem \cite[Theorem
4.2.4]{2} on rational points near planar curves by a bound which
is essentially best possible. Establishing Conjecture 2 and
Conjecture H  completes the general metric theory (i.e. the
Lebesgue and Hausdorff measure theories) for  planar curves.

 More
precisely, let $\eta<\xi$, $I=[\eta,\xi]$ and
$f:I\rightarrow{\mathbb R}$ be such that $f''$ is continuous on
$I$ and and bounded away from $0$. For convenience we suppose that
at the end points of $I$ the appropriate one sided first and
second derivatives exist. Let $\psi:\R^+\to\R^+$ be an {\it
approximating function}, that is, a real, positive decreasing
function with $\psi(t)\rightarrow0$ as $t\rightarrow\infty$, and
define, as in \cite{1},
\begin{equation}
\label{e:one1}
N_f(Q,\psi,I) \ := \ \mathrm{card} \{{\mathbf
p}/q\in{\mathbb Q}^2 \, : \, q\le Q,\, p_1/q\in I,\, |
f(p_1/q)-p_2/q|<\psi(Q)/Q\}.
\end{equation}
Here ${\mathbf p}/q := (p_1/q,p_2/q)$ with ${\mathbf p}=(p_1,p_2)
\in {\mathbb Z}^2 $ and $ q \in {\mathbb N}$.  In short, the
function $N_f(Q,\psi,I)$ counts the number of rational points with
bounded denominator lying within a specified neighbourhood of the
curve parameterized by $f$; namely  $\mathcal{C}_f := \{ (x, f(x))
\in {\mathbb R}^2 : x \in I \} \ . $
 Then firstly we show that
\begin{equation}
\label{e:one2}
N_f(Q,\psi,I) \ll \psi(Q)Q^2
\end{equation}
when $\psi(Q)\ge Q^{-\phi}$ and $\phi$ is any real number with
$0\le \phi \leq \frac23$  -- see \S\ref{cr}.  Secondly with a
further mild condition on $f$ we show that the above holds when
$\phi<1$.

Conjecture 1 of  \cite{1}, states that (\ref{e:one2}) holds for
any $f \in C^{(3)}(I)$ and  any approximating function $\psi$ such
that $t \psi(t) \to \infty$ as $ t \to \infty$. Essentially, for
$f \in  C^{(2)}(I)$ our first counting  result  requires that
$t^{2/3} \psi(t) \to \infty$ as $ t \to \infty$ and clearly falls
well short of establishing the conjecture. Nevertheless, the
result is more than adequate for establishing the stronger
$C^{(2)}$ form of Conjecture 2 of \cite{1} which states that any
$C^{(3)}$ non-degenerate planar curve is of Khinchin type for
convergence -- see \S\ref{kr}. On the other hand,  our second
counting result just falls short of establishing  Conjecture 1 in
that it essentially requires that $t^{1-\varepsilon} \psi(t) \to
\infty$ as $ t \to \infty$.  However, it is strong enough to
verify Conjecture H of \cite{1} -- the Hausdorff measure analogue
of Conjecture 2 -- see \S\ref{hh}.

\subsection{The counting results}
\label{cr}
 Let  $\eta$, $\xi$ and $f$ be as above. Furthermore, let  $\delta > 0 $
 and consider the counting function
\begin{equation}
\label{e:one4} N(Q,\delta) := {\rm{card}}\{(a,q) \in {\mathbf Z}
\times {\mathbb N} \, : \, q\le Q, \eta q<a\le \xi q,
\|qf(a/q)\|<\delta\} \ ,
\end{equation}
where $\|x\|=\min\{|x-m|:m\in\Z\}$. The main results of this paper
are

\begin{thm}
\label{t:thm1}
Suppose that $Q\ge 1$ and $0<\delta<\frac12$.  Then
\begin{equation*}
N(Q,\delta) \ll \delta Q^2 + \delta^{-\frac12}Q  \ .
\end{equation*}
\end{thm}

From this the next theorem is an easy deduction.

\begin{thm}
\label{t:thm2} Suppose that $\psi$ is an approximating function
with $\psi(Q)\ge Q^{-\phi}$ where $\phi$ is any real number with
$\phi\le\frac23$.  Then (\ref{e:one2}) holds.
\end{thm}

With a natural additional condition on $f$ we are able to extend
the validity of the bound in Theorem \ref{t:thm1}.

\begin{thm}
\label{t:thm4} Suppose that $0<\theta<1$ and
$f''\in{\rm{Lip}}_{\theta}([\eta,\xi])$ and that $Q\ge 1$ and
$0<\delta<\frac12$.  Then
\begin{equation*}
N(Q,\delta) \ll \delta Q^2 +
\delta^{-\frac12}Q^{\frac12+\varepsilon} +
\delta^{\frac{\theta-1}2}Q^{\frac{3-\theta}2}
\end{equation*}
\end{thm}

\medskip

When $\theta=1$ the proof gives the above theorem with the term
$\delta^{\frac{\theta-1}2}Q^{\frac{3-\theta}2}$ replaced by $Q\log(Q/\delta)$,
and this is then always bounded by one of the other two terms.

 We remark in passing that when $\delta>Q^{\varepsilon-1}$ our arguments can
 be extended to show that
\begin{equation*}
N(Q,\delta) \sim (\xi-\eta)\delta Q^2
\end{equation*}
and this has relevance to the further development of the Khinchin theory.
We intend to return to this in a future publication.

From Theorem \ref{t:thm4},  the next theorem is an easy deduction.

\begin{thm}
\label{t:thm5} Suppose that $0<\theta<1$ and
$f''\in{\rm{Lip}}_{\theta}([\eta,\xi])$, and suppose that $\psi$
is an approximating function with $\psi(Q)\ge Q^{-\phi}$ where
$\phi$ is any real number with $\phi \le
\frac{1+\theta}{3-\theta}$. Then (\ref{e:one2}) holds.
\end{thm}

The following statement follows immediately from Theorem
\ref{t:thm5} and essentially verifies Conjecture 1 of \cite{1}.

\begin{cor}
\label{cor1} Suppose that $f \in C^{(3)}([\eta,\xi])$, and suppose
that $\psi$ is an approximating function with $\psi(Q)\ge
Q^{-\phi}$ where $\phi$ is any real number with $\phi<1$.  Then
(\ref{e:one2}) holds.
\end{cor}

For approximating functions $\psi$ satisfying $t^{2/3} \psi(t) \to
\infty$ as $ t \to \infty$, Theorem \ref{t:thm2} removes the
factor $\delta^{-\varepsilon}$  from Huxley's estimate (see
\cite[\S1.4]{1} and \cite[Theorem 4.2.4, (4.2.20)]{2}).   With its
slightly stronger hypothesis Theorem \ref{t:thm5} also does this
for approximating functions $\psi$ satisfying
$t^{1-\varepsilon}\psi(t) \to \infty$ as $t \to \infty$ and
complements the lower bound estimate obtained in \cite[Theorem
6]{1}. Although apparently negligible, the extra factor
$\delta^{-\varepsilon}$ in Huxley's estimate renders it inadequate
for our purposes as it stands.  However, it plays an important
r\^ole in our proof. Moreover the duality described on page 72 of
Huxley \cite{2} is central to our argument. In Huxley's work the
duality occurs in an elementary way.  Here it arises as a
consequence of the harmonic analysis, where it explicitly reverses
the r\^oles of $\delta$ and $Q$.

\subsection{The Khinchin theory}

\label{kr}

Given an approximating function $\psi$, a point $\vv
y=(y_1,y_2)\in\R^2$ is called {\it simultaneously
$\psi$--approximable} if there are infinitely many $q\in {\mathbb
N}$ such that $$ \max_{1\le i\le 2}\|q y_i\|<\psi(q)  \ . $$
 Let
${\mathcal S}(\psi)$ denote the set of simultaneously
$\psi$--approximable points in $\R^2$.  Khinchin's theorem
provides a simple criteria for the `size' of ${\mathcal S}(\psi)$
expressed in terms of  two-dimensional Lebesgue measure $|\ \
|_{\R^2}$; namely
$$|{\mathcal S}(\psi)|_{\R^2} =\left\{\begin{array}{ll} \mbox{\rm
Z{\scriptsize ERO}} & {\rm if} \;\;\; \sum \;  \psi(t)^2 \;\;
<\infty\\ &
\\ \mbox{\rm F{\scriptsize ULL}} & {\rm if} \;\;\; \sum  \;  \psi(t)^2 \;\;
 =\infty \; \;
\end{array}\right.,$$
\noindent where `full' simply means that the complement of the set
under consideration is of zero measure. Now let ${\mathcal C}$ be
a planar curve and consider the set $ {\mathcal C} \cap {\mathcal
S}(\psi)   $ consisting  of points $\vv y $ on ${\mathcal C}$
which are simultaneously $\psi$--approximable. The goal is to
obtain an analogue of Khinchin's theorem for $ {\mathcal C} \cap
{\mathcal S}(\psi)   $.  Trivially, $ |{\mathcal C} \cap {\mathcal
S}(\psi)|_{\R^2} = 0  $ irrespective of the approximating function
$\psi$. Thus, when referring to the Lebesgue measure of the set $
{\mathcal C} \cap {\mathcal S}(\psi) $ it is always  with
reference to the induced Lebesgue measure $|\ \ |_{{\mathcal C}}$
on ${\mathcal C}$.  Now some useful terminology:
\begin{enumerate}
\item ${\mathcal C}$ is  of {\it Khinchin type for convergence}\/ when
$|{\mathcal C}\cap {\mathcal S}(\psi)|_{{\mathcal C}}= \mbox{{\rm
Z{\scriptsize ERO}}}$ for any approximating function $\psi$ with
$\sum\psi(t)^2<\infty$.
\\
\item ${\mathcal C}$ is  of {\it Khinchin type for divergence}\/ when
$|{\mathcal C}\cap {\mathcal S}(\psi)|_{{\mathcal C}}= \mbox{{\rm
F{\scriptsize ULL}}}$ for any approximating function $\psi$ with
$\sum\psi(t)^2=\infty$.
\end{enumerate}

\noindent To make any reasonable progress with developing a
Khinchin theory for planar curves $\mathcal{C}$, it is  reasonable
to assume that the set of points on $ \mathcal{C}$ at which the
curvature vanishes is a set of one-dimensional Lebesgue measure
zero, i.e. the curve is {\em non-degenerate }. In \cite{1}, the
following result is established.

\noindent{\bf Theorem.} {\em Any $C^{(3)}$ non-degenerate planar
curve is of Khinchin type for divergence.}

To complete the Khinchin theory for $C^{(3)}$ non--degenerate
planar curves we need to show that any such curve is of Khinchin
type for convergence. A  consequence of Theorem \ref{t:thm1}, or
equivalently a slight variant of Theorem \ref{t:thm2},  is

\begin{thm}
\label{t:thm3} Any $C^{(2)}$ non-degenerate planar curve
 is of Khinchin type for convergence.
\end{thm}

In the case $\psi:t\to t^{-v}$ with $v>0$, let us write ${\mathcal
S}(v)$ for ${\mathcal S}(\psi)$. Note that in view of Dirichlet's
theorem (simultaneous version), ${\mathcal S}(v)= \R^2 $ for any
$v \leq 1/2$ and so $|\mathcal{C} \cap
\mathcal{S}(v)|_{\mathcal{C}} = |\mathcal{C}|_{\mathcal{C}} :=
\mbox{F{\scriptsize ULL}} $ for any $v \leq 1/2$. It is easily
verified that Theorem \ref{t:thm3} implies the following
`extremality' result due to  Schmidt \cite{sch}.

\noindent{\bf Corollary (Schmidt).} {\em Let ${\mathcal C}$ be a
$C^{(2)}$ non-degenerate planar curve. Then, for  any $ v
> 1/2$ $$ |\mathcal{C}
\cap \mathcal{S}(v)|_{\mathcal{C}} \ = \ 0  
 \ . $$ }

\noindent To be precise, Schmidt actually requires that ${\mathcal
C} $ is a  $C^{(3)}$ non-degenerate planar  curve.   
For further background, including a
comprehensive account of related works,  we refer the reader to
\cite[\S1]{1}.

\subsection{The Jarn\'{\i}k theory}
\label{hh}

Jarn\'{\i}k's theorem is a Hausdorff measure version of Khinchin's
theorem in that it  provides a simple criteria for the `size' of
${\mathcal S}(\psi)$ expressed in terms of  $s$--dimensional
Hausdorff measure $\hs$.
 The Hausdorff measure and  dimension of a set $X \in \R^2$ is
defined as follows. For $\rho > 0$, a countable collection $
\left\{B_{i} \right\} $ of Euclidean balls  in  $\R^{2}$ with
diameter ${\rm diam}(B_i)  \leq \rho $ for each $i$  such that $X
\subset \bigcup_{i} B_{i} $ is called a $ \rho $-cover for $X$.
Let $s$ be a non-negative number and define $
 \hs_{ \rho } (X)
  \; = \; \inf \left\{ \textstyle{\sum_{i}} {\rm diam}(B_i)^s
\ :   \{ B_{i} \}  {\rm \  is\ a\  } \rho {\rm -cover\  of\ } X
\right\}  $, where the infimum is taken over all possible $ \rho
$-covers of $X$. The $s$-{\it dimensional Hausdorff measure}
$\hs(X)$ is defined by $$ \hs(X) := \lim_{\rho \to 0}
\hs_{\rho}(X) = \sup_{\rho > 0} \hs_{\rho}(X) \  $$ and the  {\it
Hausdorff dimension} $\dim X$ of $X$ is defined by $$ \dim \, X :=
\inf \left\{ s : \hs (X) =0 \right\} = \sup \left\{ s : \hs (X) =
\infty \right\} \, . $$

Jarn\'{\i}k's theorem shows that the  $s$--dimensional Hausdorff
measure $\hs({\mathcal S}(\psi))$  of the set ${\mathcal S}(\psi)$
satisfies an elegant `zero-infinity' law. Let $s \in (0,2)$ and
$\psi$ be  an approximating function. Then $$ \hs\left({\mathcal
S}(\psi)\right)=\left\{\begin{array}{ll} 0 & {\rm when} \;\;\;
\sum \; t^{2-s} \,  \psi(t)^s \;\;
 <\infty\\ &
\\ \infty & {\rm when} \;\;\; \sum \;
 t^{2-s} \,  \psi(t)^s  \;\;  =\infty
\end{array}\right..$$

\noindent Note that this  trivially  implies that  $ \dim
{\mathcal S}(\psi)  =  \inf \{ s : \mbox{$\sum $} \; t^{2-s} \,
\psi(t)^s < \infty \} $.

\medskip

Now let ${\mathcal C}$ be a planar curve.  The goal is to obtain
an analogue of Jarn\'{\i}k's theorem for $ {\mathcal C} \cap
{\mathcal S}(\psi)   $.  In particular, our aim is to establish the
following conjecture stated in \cite{1}.

\noindent{{\bf Conjecture  H }} {\em Let $s \in (1/2,1)$ and
$\psi$ be an approximating function. Let $f\in C^{(3)}(I)$, where
$I$ is an interval and let ${\mathcal C}_f : = \{(x,f(x)):x\in I
\} $. Assume that $ \dim\{x\in I: f''(x)=0\} \leq 1/2 $.  Then
$$
\hs\left({\mathcal C}_f\cap {\mathcal
S}(\psi)
\right)=\left\{
\begin{array}{ll} 0 & {\rm when} \;\;\; \sum
\; t^{1-s} \,  \psi(t)^{s+1} \;\;
 <\infty\\ &
\\ \infty & {\rm when} \;\;\; \sum \;
 t^{1-s} \,  \psi(t)^{s+1}  \;\;  =\infty
\end{array}
\right..$$  }

\smallskip

The divergent part of the above statement, namely that
$$\hs\left({\mathcal C}_f\cap {\mathcal
S}(\psi) \right)= \infty \quad \text{ when } \quad \sum \;
 t^{1-s} \,  \psi(t)^{s+1}  \;\;  =\infty,$$
is Theorem 3 in \cite{1},
and  so the main substance of the conjecture is the
convergence part. A consequence of Theorem \ref{t:thm4} above, or
equivalently a slight variant of Corollary \ref{cor1}, is the
completion of the proof of Conjecture H.

\begin{thm}
\label{t:thm6} Let $s \in (1/2,1)$ and
$\psi$ be an approximating function. Let $f\in C^{(3)}(I)$, where
$I$ is an interval and let ${\mathcal C}_f : = \{(x,f(x)):x\in I
\} $. Assume that $ \dim\{x\in I: f''(x)=0\} \leq 1/2 $.  Then
$$ \hs\left({\mathcal C}_f\cap {\mathcal S}(\psi)\right) \ = \ 0
\quad {\rm when} \quad  \textstyle{ \sum } \; t^{1-s} \,
\psi(t)^{s+1}
 <\infty \ .$$
\end{thm}

For further background, including an explanation of the conditions
in Conjecture H and a  comprehensive account of related works, we
refer the reader to \cite[\S1]{1}.

\Section{The proof of Theorem \ref{t:thm1}}

\noindent It clearly suffices to prove Theorem \ref{t:thm1} and
indeed Theorem \ref{t:thm4} with $N(Q,\delta)$ replaced by $$
\widetilde{N}(Q,\delta) := {\rm{card}}\{(a,q) \in {\mathbf Z}
\times {\mathbb N} \, : \, Q < q\le 2Q, \eta q<a\le \xi q,
\|qf(a/q)\|<\delta\} \ . $$ 

 \noindent Let
\begin{equation}
\label{e:two1}
J=\left\lfloor
\frac1{2\delta}
\right\rfloor
\end{equation}
and consider the Fej\'er kernel
\begin{equation*}
\mathcal K_J(\alpha) = J^{-2}\left|
\sum_{h=1}^Je(h\alpha)
\right|^2  = \left(
\frac{\sin\pi J\alpha}{J\sin\pi\alpha}
\right)^2.
\end{equation*}
When $\|\alpha\|\le \delta$ we have $|\sin\pi J\alpha| =
\sin\pi\|J\alpha\|\ge 2\| J\alpha\| = 2\|\,
J\|\alpha\|\, \|=2J\|\alpha\|$, since $J\|\alpha\| \le \delta
\left\lfloor \frac1{2\delta} \right\rfloor \le \frac12$.  Hence,
when $\|\alpha\|\le\delta$, we have
\begin{equation*}
\mathcal K_J(\alpha) \ge \frac{2\|\alpha\|J}{J\pi\|\alpha\|}=\frac2{\pi}.
\end{equation*}
Thus
\begin{equation*}
\widetilde{N}(Q,\delta)\le \frac{\pi}2 \sum_{Q<q\le 2Q} \sum_{\eta
q < a \le \xi q} \mathcal K_J \big( qf(a/q) \big).
\end{equation*}
Since
\begin{equation*}
\mathcal K_J(\alpha) = \sum_{j=-J}^J \frac{J-|j|}{J^2} \,
e(j\alpha)
\end{equation*}
we have
\begin{equation*}
\widetilde{N}(Q,\delta)\le \pi\delta (\xi-\eta)Q^2 + N_1 +O(\delta
Q) = N_1 +O(\delta Q^2)
\end{equation*}
where $$N_1= \frac{\pi}2 \sum_{0<|j|\le J} \frac{J-|j|}{J^2}
\sum_{Q<q\le 2Q} \  \sum_{\eta q<a\le \xi q}
e\big(jqf(a/q)\big).$$ We observe that the function $F(\alpha)=
jqf(\alpha/q)$ has derivative $jf'(\alpha/q)$. Given $j$ with
$0<|j|\le J$ we define $$H_-=\lfloor\inf jf'(\beta)\rfloor-1,\ \
H_+=\lceil\sup jf'(\beta)\rceil+1,$$ $$h_-=\lceil\inf
jf'(\beta)\rceil+1,\ \ h_+=\lfloor\sup jf'(\beta)\rfloor-1$$ where
the extrema are over the interval $[\eta,\xi]$.  Then, by Lemma
4.2 of Vaughan \cite{4}, $$\sum_{\eta q<a\le \xi q}
e\big(jqf(a/q)\big) = \sum_{H_-\le h\le H_+} \int_{\eta q}^{\xi q}
e\big(jqf(\alpha/q)-h\alpha\big)d\alpha +O\big(\log(2+H)\big)$$
where $H=\max(|H_-|,|H_+|)$.  Clearly $H\ll|j|\le J$ and so
$$N_1=N_2 + O\big(Q\log\textstyle{\frac1\delta}\big)$$ where $$N_2
= \frac{\pi}2 \sum_{0<|j|\le J} \frac{J-|j|}{J^2} \sum_{Q<q\le 2Q}
\sum_{H_-\le h\le H_+} \int_{\eta q}^{\xi q}
e\big(jqf(\alpha/q)-h\alpha\big)d\alpha.$$ The integral here is
$$q\int_{\eta}^{\xi} e\big(q(jf(\beta)-h\beta)\big)d\beta.$$ The
function $g(\beta)= q(jf(\beta)-h\beta)$ has second derivative
$qjf''(\beta)$ whose modulus lies between constant multiples of
$q|j|$.  Hence, by Lemma 4.4 of Titchmarsh \cite{3}, for any
subinterval $\mathcal I$ of $[\eta,\xi]$,
\begin{equation}
\label{e:two2}
\int_{\mathcal I} e\big(q(jf(\beta)-h\beta)\big)d\beta \ll \frac1{\sqrt{q|j|}}.
\end{equation}
Thus the contribution to $N_2$ from any $h$ with $H_-\le h\le h_-$ or
$h_+\le h \le H_+$ is
$$\ll J^{-1} \sum_{j=1}^J j^{-1/2} \sum_{Q<q\le 2Q} q^{1/2}.$$
Therefore
$$N_2=N_3+O\left(
\delta^{\frac12}Q^{\frac32}
\right).$$
where
\begin{equation}
\label{e:two3}
N_3= \frac{\pi}2 \sum_{0<|j|\le J} \frac{J-|j|}{J^2} \sum_{Q<q\le 2Q}
q \sum_{h_-< h < h_+} \int_{\eta}^{\xi} e\big(q(jf(\beta)-h\beta)\big)d\beta.
\end{equation}
The sum over $h$ here is taken to be empty when $h_+\le h_-+1$.
\par
We have $$\delta^{\frac12}Q^{\frac32}=\big(\delta
Q^2\big)^\frac12\big(Q\big)^\frac12\le \delta Q^2 +Q.$$ Thus it
remains to treat $N_3$.

Since $f'$ is continuous and $\inf jf'(\beta) <h_-< h <h_+< \sup
jf'(\beta)$ it follows that there is a
$\beta_h=\beta_{j,h}\in[\eta,\xi]$ such that $jf'(\beta_h)=h$.
Let
\begin{equation}
\label{e:two4}
\lambda_h=\lambda_{j,h}=\|jf(\beta_h)-h\beta_h\|
\end{equation}
We need to bound various sums involving $\lambda_h$.  To that end
the following lemma is very useful.
\begin{lem}
\label{l:lem1} Suppose that $\phi$ has a continuous second
derivative on $[\Upsilon,\Xi]$ which is bounded away from $0$, and
suppose that $\Psi$ is real and satisfies $0< \Psi< \frac14$. Then
for any fixed $\varepsilon>0$ and $R\ge 1$, the number $M$ of
triples of integers $r,b,c$ such that $(r,b,c)=1$, $R\le r< 2R$,
$\Upsilon r<b\le \Xi r$ and  $|r\phi(b/r)-c|\le \Psi$ satisfies
$$M\ll_{\varepsilon} \Psi^{1-\varepsilon} R^2+R.$$
\end{lem}
\begin{proof}
If $\Upsilon<0<\Xi$, then we split $[\Upsilon,\Xi]$ into two
subintervals $[\Upsilon,0]$, $[0,\Xi]$ and consider them
separately.  Thus we may suppose $0\notin(\Upsilon,\Xi)$.  If
$\Xi\le 0$, then by replacing $b/r$ by $-b/r$ and $\Psi(\alpha)$
by $\Psi(-\alpha)$ we can transfer our attention to the interval
$[-\Xi,-\Upsilon]$.  Thus it always suffices to consider intervals
$[\Upsilon,\Xi]$ with $0\le\Upsilon\le\Xi$.  Now choose
$K\in\mathbb N$ so that $K>\Xi$, say $K=\lfloor \Xi\rfloor +1$. We
extend the definition of $\phi$ so that $\phi$ is twice
differentiable with a continuous second derivative and  bounded
away from $0$ on the whole of $[0,1]$.  For example, if
$\Upsilon/K>0$, then for $0\le \alpha<\Upsilon/K$ we can take
$\phi(\alpha)= \frac12(\alpha-\Upsilon)^2\phi''(\Upsilon)
+(\alpha-\Upsilon)\phi'(\Upsilon) +  \phi(\Upsilon)$, and likewise
when $\Xi<\alpha\le K$.  If we now define $F(x)$ on $[0,1]$ by
$F(x)=\phi(xK)/K$, then $F$ will satisfy the hypothesis of Theorem
4.2.4 of Huxley \cite{2}.  The condition $(r,b,c)=1$ ensures that
the rational number points $(b/r,c/r)$ are counted uniquely.  Thus
$M$ is bounded by the number of $r$, $b$, $c$ with $R\le r< 2R$,
$0<b<rK$, $|\phi(b/r)-c/r|\le R^{-1}\Psi$.  We take $T=K$, $M=K$,
$\Delta=K^{-1}$, $Q=R$, $\delta=\Psi$ and apply the conclusion
(4.2.20), {\it ibidem}, to obtain the desired result.
\end{proof}
\begin{lem}
\label{l:lem2} Suppose that $\phi$ has a continuous second
derivative on $[\Upsilon,\Xi]$ which is bounded away from $0$, and
suppose that $\Psi$ is real and satisfies $0< \Psi< \frac14$. Then
for any $\varepsilon>0$ and $R\ge 1$, $$\sum_{R\le r< 2R}  \ \
\sum_{\substack{\Upsilon r<b\le \Xi r \\ \|r\phi(b/r)\|\le \Psi}}
1 \ll_{\varepsilon} \Psi^{1-\varepsilon} R^2+R\log 2R.$$
\end{lem}
\begin{proof}
For a given pair $r$, $b$ counted in the double sum let $c$ be the
unique integer with $|r\phi(b/r)-c|\le \Psi$.  We sort the triples
according to the value of the greatest common divisor $(r,b,c)=d$,
say.  Then the double sum does not exceed $$\sum_{d\le 2R}M(d)$$
where $M(d)$ is the number of triples of integers $s,g,h$ such
that $(s,g,h)=1$, $R/d\le s< 2R/d$, $\Upsilon s<g\le \Xi s$ and
$|s\phi(g/s)-h|\le \Psi d^{-1}$.  By Lemma 2.1,
$$M(d)\ll_{\varepsilon} \big(\Psi d^{-1})^{1-\varepsilon}(R/d)^2 +
R/d.$$ Summing over the $d\le 2R$ gives the lemma.
\end{proof}

We apply this through the next lemma.
\begin{lem}
\label{l:lem3}
We have
\begin{equation}
\label{e:two5} \sum_{0<|j|\le J} \  \sum_{\substack{h_-<h<h_+ \\
\lambda_h> Q^{-1}}} |j|^{-\frac12} \lambda_h^{-\frac12} \ll
J^{\frac32}+J^{\frac12}(\log J)Q^{\frac12},
\end{equation}
\begin{equation}
\label{e:two6} \sum_{0<|j|\le J} \ \sum_{\substack{h_-<h<h_+ \\
\lambda_h> Q^{-1}}} |j|^{-\frac12} \lambda_h^{-1} \ll
J^{\frac32}Q^{\varepsilon}+J^{\frac12}(\log J)Q
\end{equation}
and
\begin{equation}
\label{e:two7} \sum_{0<|j|\le J} \  \sum_{\substack{h_-<h<h_+ \\
\lambda_h\le Q^{-1}}} |j|^{-\frac12}  \ll
J^{\frac32}Q^{\varepsilon-1}+J^{\frac12}\log J.
\end{equation}
\end{lem}
\begin{proof}
Clearly in (\ref{e:two5}) and (\ref{e:two6}) we can restrict our
attention to terms with $\lambda_h<\frac14$, since those terms
with $\lambda_h\ge \frac14$ contribute $\ll J^{\frac32}$ to the
total. Let $$\Upsilon=\inf f'(\beta),\,\Xi=\sup f'(\beta)$$ where
the extrema are taken over $[\eta,\xi]$.  When $j<0$ we replace
$j$ by $-j$ and $h$ by $-h$ in each of the sums in question, and
write $\beta_h$ for $\beta_{-h}$ and $\lambda_h$ for
$\lambda_{-h}$ to see that the sums are bounded by $$\sum_{j=1}^J
 \ \sum_{\substack{\Upsilon j<h<\Xi j \\ \lambda_h> Q^{-1}}}
j^{-\frac12} \lambda_h^{-\frac12},$$ $$\sum_{j=1}^J \
\sum_{\substack{\Upsilon j<h<\Xi j \\ \lambda_h> Q^{-1}}}
j^{-\frac12} \lambda_h^{-1},$$ and $$\sum_{j=1}^J \
\sum_{\substack{\Upsilon j<h<\Xi j \\ \lambda_h\le Q^{-1}}}
j^{-\frac12}$$ respectively.

Let $g$ denote the inverse function of $f'$, so that $g$ is
defined on $[\Upsilon,\Xi]$ and $\beta_h=g(h/j)$.  Let $F(\alpha)=
\alpha g(\alpha)-f\big(g(\alpha)\big)$.  Then $$F'(\alpha)= \alpha
g'(\alpha) + g(\alpha) - f'\big(g(\alpha)\big)g'(\alpha) =
g(\alpha)$$ and $$F''(\alpha) =
g'(\alpha)=\frac1{f''(g(\alpha))}$$ and so, in particular, $F''$
is bounded away from $0$.  Thus $$\lambda_h=\|jF(h/j)\|$$ and $F$
satisfies the conditions on $\phi$ in Lemma 2.2.  The desired
bounds now follow by partial summation.
\end{proof}

We now return to the estimation of $N_3$, defined by
(\ref{e:two3}).  By (\ref{e:two2}), the terms in $N_3$ with
$\lambda_h\le Q^{-1}$ contribute $$\ll J^{-1}\sum_{0<|j|\le J}
 \ \sum_{\substack{h_-<h<h_+ \\ \lambda_h\le Q^{-1}}}  \ \sum_{Q<q\le
2Q} q^{\frac12}|j|^{-\frac12}$$ and by (\ref{e:two7}) this is
$$\ll J^{\frac12}Q^{\frac12+\varepsilon} + J^{-\frac12}(\log
J)Q^{\frac32}.$$ Hence
\begin{equation}
\label{e:two8} N_3=N_4 +O\big(\delta^{-\frac12}
Q^{\frac12+\varepsilon} +
\delta^{\frac12}\big(\log\textstyle{\frac1\delta}\big)Q^{\frac32}\big)
\end{equation}
where
$$N_4=\frac{\pi}2 \sum_{0<|j|\le J} \frac{J-|j|}{J^2}
\sum_{Q<q\le 2Q} q  \ \sum_{\substack{h_-< h < h_+ \\
\lambda_h>Q^{-1}}} \int_{\eta}^{\xi}
e\big(q(jf(\beta)-h\beta)\big)d\beta.$$
Since
\begin{equation}
\label{e:two9}
\delta^{\frac12}\left(\log\textstyle{\frac1\delta}\right)Q^{\frac32}=\big(\delta
Q^2\big)^\frac12 \left(\left(\log\textstyle{\frac1\delta}\right)^2
Q\right)^\frac12 \le \delta Q^2 +\delta^{-\frac12}Q
\end{equation}
this gives
$$N_3=N_4+O\big(\delta Q^2+\delta^{-\frac12}Q\big).$$
\par
Let $c=\big(\sup{|f''(\beta)|}\big)^{-1/2}$ where the supremum is
taken over $[\eta,\xi]$. The set $\mathcal A(j,h)$ of those
$\beta$ in $[\eta,\xi]$ for which
$|\beta-\beta_h|>c\sqrt{\lambda_h/|j|}$ consists of at most two
intervals, and may be empty. By the mean value theorem, for such
$\beta$ we have $$jf'(\beta)-h = (\beta-\beta_h)jf''(\beta^*)$$
for some $\beta^*\in[\eta,\xi]$.  Thus $$|jf'(\beta)-h|\gg
\sqrt{|j|\lambda_h}.$$ Hence, by integration by parts, we have
$$\int_{\mathcal A(j,h)} e\big(q(jf(\beta)-h\beta)\big)d\beta \ll
\frac1{q\sqrt{|j|\lambda_h}}.$$ Therefore the total contribution
to $N_4$ from the $\mathcal A(j,h)$ is $$\ll J^{-1}Q
\sum_{0<|j|\le J} \sum_{\substack{h_-<h<h_+ \\ \lambda_h> Q^{-1}}}
\frac1{\sqrt{|j|\lambda_h}}$$ and by (\ref{e:two5}) this is $$\ll
J^{\frac12}Q + J^{-\frac12}(\log J)Q^{\frac32}.$$ Thus, by
(\ref{e:two9}), $$N_4=N_5+O\big(\delta Q^2 +
\delta^{-\frac12}Q\big)$$ where $$N_5=\frac{\pi}2 \sum_{0<|j|\le
J} \frac{J-|j|}{J^2} \sum_{Q<q\le 2Q} q \sum_{\substack{h_-< h <
h_+ \\ \lambda_h>Q^{-1}}} \int_{\mathcal B(j,h)}
e\big(q(jf(\beta)-h\beta)\big)d\beta$$ and $\mathcal B(j,h)$
denotes the set of $\beta\in[\eta,\xi]$ with $|\beta-\beta_h|\le
c\sqrt{\lambda_h/|j|}$.

Given $j$ and $h$ included in the sum, choose $n=n(j,h)$ so that
$\lambda_h=|jf(\beta_h)-h\beta_h-n|$.  For $\beta\in B(j,h)$ we
have
\begin{equation}
\label{e:two10} jf(\beta)-h\beta-n=jf(\beta_h)-h\beta_h-n +
\textstyle{\frac12} (\beta-\beta_h)^2 j f''(\beta^\flat)
\end{equation}
where $\beta^\flat\in[\eta,\xi]$.  When $\frac14\le\lambda_h$ we
have $$\textstyle{\frac18}\le \textstyle{\frac12} \lambda_h \le
|jf(\beta) - h\beta-n| \le \textstyle{\frac32} \lambda_h\le
\textstyle{\frac34}.$$ Thus $\|jf(\beta)-h\beta\| =
|jf(\beta)-h\beta-m|$ with $m=n$ or $m=n\pm1$, and so $$\frac18\le
\|jf(\beta)-h\beta\|.$$ On the other hand, when
$\lambda_h<\frac14$ the identity (\ref{e:two10}) shows that
$$\textstyle\frac12\lambda_h \le \|jf(\beta)-h\beta\| \le
\textstyle\frac32\lambda_h$$ and so generally
$$\|jf(\beta)-h\beta\|\asymp \lambda_h.$$ Therefore for $j$ and
$h$ included in the sum we have $$ \int_{\mathcal
B(j,h)}\sum_{Q<q\le 2Q} q e\big(q(jf(\beta)-h\beta)\big)d\beta \ll
Q\lambda_h^{-1}\mathrm{meas}\mathcal B(j,h) \ll  Q\lambda_h^{-1}
\sqrt{\lambda_h/|j|}. $$ and hence $$N_5\ll J^{-1} Q\sum_{0<|j|\le
J} \sum_{\substack{h_-< h < h_+ \\ \lambda_h>Q^{-1}}}
|j|^{-\frac12}\lambda_h^{-\frac12}. $$ Thus, by (\ref{e:two5}),
$$N_5\ll J^{-1}Q\big(J^{\frac32} + Q^{\frac12}J^{\frac12}\log
J\big) \ll \delta^{-\frac12}Q +
\delta^{\frac12}\big(\log\textstyle{\frac1\delta}\big)Q^{\frac32}.$$

\noindent This with (\ref{e:two9}) completes the proof of Theorem
\ref{t:thm1}.

\Section{ The proof of Theorem \ref{t:thm2}}
\noindent By (\ref{e:one1}), when $\psi(Q)\le \frac12$, $N_f(Q,\psi,I)$ is
$$\le \mathrm{card}\{a,q \, : \, q\le Q,\, a\in qI,\,
\|qf(a/q)\|<q\psi(Q)/Q\}$$
and this is bounded by
$$\mathrm{card}\{a,q \, : \, q\le Q,\, a\in qI,\,
\|qf(a/q)\|< \psi(Q)\} \  . $$ Now the conclusion is immediate
from Theorem \ref{t:thm1}.

\Section{The proof of Theorem \ref{t:thm4}}

\noindent For convenience we extend the definition of $f$ to
${\mathbb R}$ by defining $f(\beta)$ to be
$\frac12(\beta-\xi)^2f''(\xi)+(\beta-\xi)f'(\xi)+f(\xi)$ when
$\beta>\xi$ and to be
$\frac12(\beta-\eta)^2f''(\eta)+(\beta-\eta)f'(\eta)+f(\eta)$ when
$\beta<\eta$.  Note that then $f''\in\rm{Lip}_{\theta}(\mathbb R)$
and $f''$ is still bounded away from $0$ and is bounded .\par We
follow the proof of Theorem \ref{t:thm1} as far as (\ref{e:two8}).  We note
that the complete error term here is in fact
$$\delta Q^2 + Q\log\textstyle{\frac1\delta} + \delta^{-\frac12}Q^{\frac12+\varepsilon}
+\delta^{\frac12}(\log\textstyle{\frac1\delta})Q^{\frac32}$$
Thus
$$\widetilde{N}(Q,\delta) \ll N_4 +
\delta Q^2 + Q\log\textstyle{\frac1\delta} + \delta^{-\frac12}Q^{\frac12+\varepsilon}
+ \delta^{\frac12}(\log\textstyle{\frac1\delta})Q^{\frac32}$$
where
$$N_4= \frac{\pi}2 \sum_{0<|j|\le J} \frac{J-|j|}{J^2}
\sum_{Q<q\le 2Q} \! q \  \sum_{\substack{h_-< h < h_+ \\
\lambda_h>Q^{-1}}} \int_{\eta}^{\xi}
e\big(q(jf(\beta)-h\beta)\big)d\beta.$$ Moreover, given $j$ and
$h$ included in the sums there is unique
$\beta_h=\beta_{j,h}$ such that $$f'(\beta_h)=h/j.$$ Let
$$\mu=\frac{\xi-\eta}2.$$ Then in the integral above we replace
the interval $[\eta,\xi]$ by $[\beta_h-\mu,\beta_h+\mu]$.  For any
$\beta$ not in both intervals we have $|\beta-\beta_h|\ge \mu$,
$\beta\le\eta$, or $\beta\ge\xi$.  For some $\beta^*\in[\eta,\xi]$
we have $(\beta_h-\eta)jf''(\beta^*)=jf'(\beta_h)-jf'(\eta)\ge h
-h_-$ so $\beta_h-\eta\gg (h-h_-)/|j|$ and likewise $\xi-\beta_h
\gg (h_+-h)/|j|$.  Hence, if $\beta\le\eta$, then
$\beta_h-\beta\gg (h-h_-)/|j|$, and if $\beta\ge\xi$, then
$\beta-\beta_h \gg (h_+-h)/|j|$. Moreover, as $\mu\gg (h-h_-)/|j|$
and $\mu\gg (h_+-h)/|j|$ it follows that whenever $\beta$ is not
in both intervals we have either $|\beta-\beta_h|\gg (h-h_-)/|j|$ or
$|\beta-\beta_h| \gg (h_+-h)/|j|$.  For any such $\beta$ there is a
$\beta^\flat$ such that $jf'(\beta)-h = j(f'(\beta)-f'(\beta_h)) =
j(\beta-\beta_h)f''(\beta^\flat)$, whence $|jf'(\beta)-h|\gg h-h_-$
or $|jf'(\beta)-h|\gg h_+-h$.  It then follows by integration by parts
that if ${\mathcal
A}=[\eta,\xi]\backslash[\beta_h-\mu,\beta_h+\mu]$ or ${\mathcal
A}=[\beta_h-\mu,\beta_h+\mu]\backslash[\eta,\xi]$, then
$$\int_{\mathcal A}e\big(q(jf(\beta)-h\beta)\big)d\beta \ll
\frac1{q(h-h_-)}+\frac1{q(h_+-h)}. $$
Thus
$$N_4=N_5+O\left(
\sum_{0<|j|\le J}
\sum_{h_-<h<h_+}\frac{Q/J}{h-h_-}+\frac{Q/J}{h_+-h} \right)$$
where
$$N_5=\frac{\pi}2 \sum_{0<|j|\le J} \frac{J-|j|}{J^2}
\sum_{\substack{h_-< h < h_+ \\ \lambda_h>Q^{-1}}} \sum_{Q<q\le
2Q} q \int_{\beta_h-\mu}^{\beta_h+\mu}
e\big(q(jf(\beta)-h\beta)\big)d\beta.$$
Thus
$$N_4=N_5+O\big(Q\log\textstyle{\frac1\delta}\big).$$
For convenience we write
$$F(\alpha)=F(\alpha;j,h)=\big(f(\alpha+\beta_h)-f(\beta_h)\big)-h\alpha/j.$$
Then $$F(0)=0,\quad F'(\alpha)=f'(\alpha+\beta_h)-h/j,\quad
F'(0)=0,\quad F''(\alpha)=f''(\alpha+\beta_h),$$
and
$$\int_{\beta_h-\mu}^{\beta_h+\mu}
e\big(q(jf(\beta)-h\beta)\big)d\beta = e(q\phi_h)\int_{-\mu}^{\mu}
e\big(qjF(\alpha)\big)d\alpha$$
where
$$\phi_h=\phi_{j,h}=jf(\beta_h)-h\beta_h$$
so that
$$\lambda_h=\|\phi_h\|.$$
Since $f''\in\rm{Lip}_\theta(\mathbb R)$, we have $F''\in\rm{Lip}_\theta(\mathbb R)$
and so, in particular,
$$F''(\alpha)=F''(0)+O\big(|\alpha|^{\theta}\big) =
f''(\beta_h)+O\big(|\alpha|^{\theta}\big),$$
and thus
$$F'(\alpha)=\alpha f''(\beta_h) +
O\big(|\alpha|^{1+\theta}\big),\quad F(\alpha)=\textstyle{\frac12}\alpha^2
f''(\beta_h) + O\big(|\alpha|^{2+\theta}\big).$$
For brevity write
$c_2=f''(\beta_h)$.
\par
Since $f''$, and hence $F''$, is bounded
and bounded away from $0$, and $f''$ is continuous it follows that
$F'$ is strictly monotonic and so can only change sign once.  But
$F'(0)=0$.  We suppose for the time being that $c_2>0$.  Now $F'$
is strictly increasing, and hence positive when $\alpha>0$.  Thus
$F$ is strictly increasing for $\alpha\ge0$ and positive for $\alpha>0$.
Let $G$ be the
inverse function of $F$ on $[0,\infty)$.  Then $G'$ exists on
$(0,\infty)$ and $G'(\beta)=1/F'\big(G(\beta)\big)$.  Thus for any
$\nu$ with $$0<\nu<\mu$$ we have $$\int_{\nu}^{\mu}
e\big(qjF(\alpha)\big)d\alpha = \int_{F(\nu)}^{F(\mu)}
e(qj\beta)G'(\beta) d\beta.$$ Note that we will eventually choose
$\nu$ to be judicially small in terms of $q$ and $j$. Since $F'$
is non-zero for $\alpha>0$ it follows that $G''$ exists on
$(0,\infty)$, and is continuous, and so by integration by parts we have
$$\int_{F(\nu)}^{F(\mu)} e(qj\beta)G'(\beta) d\beta = \left[
\frac{e(qj\beta)G'(\beta)}{2\pi iqj} \right]_{F(\nu)}^{F(\mu)}
-\int_{F(\nu)}^{F(\mu)} \frac{e(qj\beta)}{2\pi iqj}G''(\beta)
d\beta.$$ Moreover
$$G''(\beta)=-\frac{F''(G(\beta))G'(\beta)}{F'(G(\beta))^2} = -
\frac{F''(G(\beta))}{G'(\beta)^3}.$$ We also have, for $\alpha>0$
$$\beta=F(\alpha)=\textstyle{\frac12}c_2\alpha^2 +
O\big(\alpha^{2+\theta}\big).$$ Since $\mu\ll1$ it follows that
for $0<\alpha\le \mu$ we have
$$G(\beta)=\alpha=\sqrt{\frac{2\beta}{c_2}}\left(
1+O\big(\beta^{\theta/2}\big) \right)= \sqrt{\frac{2\beta}{c_2}} +
O\big(\beta^{(1+\theta)/2}\big).$$ We further have
$$F'\big(G(\beta)\big) =F'(\alpha) = \sqrt{2c_2\beta} +
O\big(\beta^{(1+\theta)/2}\big).$$ and $$F''\big(G(\beta)\big) =
c_2+O\big(\alpha^{\theta}\big) = c_2 +
O\big(\beta^{\theta/2}\big).$$ Hence $$G'(\beta) =
\frac1{F'\big(G(\beta)\big)} = (2c_2\beta)^{-\frac12} +
O\big(\beta^{(\theta-1)/2}\big)$$ and $$G''(\beta) = - \left(
c_2+O\big(\beta^{\theta/2}\big) \right)\left(
\big(2c_2\beta\big)^{-\frac12}+O\big(\beta^{(\theta-1)/2}\big)
\right)^3 = -\frac{c_2}{(2c_2\beta)^{3/2}} +
O\big(\beta^{(\theta-3)/2}\big).$$
\par
Substituting the above approximations we have
\begin{align*}
&\left[ \frac{e(qj\beta)G'(\beta)}{2\pi iqj}
\right]_{F(\nu)}^{F(\mu)} - \ \int_{F(\nu)}^{F(\mu)}
\frac{e(qj\beta)}{2\pi iqj}G''(\beta) d\beta\\ \\ &= -
\frac{e(qjF(\nu))G'(F(\nu))}{2\pi iqj} + \int_{F(\nu)}^{\infty}
\frac{e(qj\beta)}{2\pi iqj}\frac{c_2}{(2c_2\beta)^{3/2}}d\beta +
E\\
\end{align*}
where
$$E\ll \frac1{q|j|} +\int_{F(\nu)}^{F(\mu)} \frac{\beta^{(\theta-3)/2}}{q|j|}d\beta \ll
\frac{F(\nu)^{(\theta-1)/2}+1}{q|j|}\ll \frac{\nu^{\theta-1}+1}{q|j|}.$$
We also have
$$G'(F(\nu))= (2c_2F(\nu))^{-\frac12} + O\big(F(\nu)^{(\theta-1)/2}\big).$$
Hence, by substitution and integration by parts,
\begin{align*}
\left[ \frac{e(qj\beta)G'(\beta)}{2\pi iqj}
\right]_{F(\nu)}^{F(\mu)} &- \ \int_{F(\nu)}^{F(\mu)}
\frac{e(qj\beta)}{2\pi iqj}G''(\beta) d\beta\\ \\  &=
\int_{F(\nu)}^{\infty} \frac{e(qj\beta)}{\sqrt{2c_2\beta}} d\beta
+ O\left( \frac{\nu^{\theta-1}+1}{q|j|} \right).\\
\end{align*}
We now turn to $$\int_0^{\nu} e\big(qjF(\alpha)\big)d\alpha.$$
This differs from $$\int_0^{\nu}
e\big(qj{\textstyle{\frac12}}c_2\alpha^2)\big)d\alpha =
\int_0^{\textstyle{\frac12}c_2\nu^2}
\frac{e(qj\beta)}{\sqrt{2c_2\beta}}d\beta$$ by $$\ll
\int_0^{\nu}q|j|\alpha^{2+\theta}d\alpha \ll q|j|\nu^{3+\theta}.$$
Now $F(\nu)= \textstyle{\frac12}c_2\nu^2 +
O\big(\nu^{2+\theta}\big)$ and so
$$\int_{\textstyle{\frac12}c_2\nu^2}^{F(\nu)}
\frac{e(qj\beta)}{\sqrt{2c_2\beta}}d\beta \ll \nu^{1+\theta}.$$
The choice $\nu=c/\sqrt{q|j|}$, where the positive constant $c$ is
chosen to ensure that $\nu<\mu$, gives $$\int_0^{\mu}
e\big(qjF(\alpha)\big)d\alpha = \int_0^{\infty}
\frac{e(qj\beta)}{\sqrt{2c_2\beta}} d\beta +
O\big((q|j|)^{(-1-\theta)/2}\big).$$ Hence $$\int_0^{\mu}
e\big(qjF(\alpha)\big)d\alpha=\frac{W_{{\rm{sgn}}(j)}}{\sqrt{qc_2|j|}}
+ O\big((q|j|)^{(-1-\theta)/2}\big)$$ where
$$W_{\pm}=\int_0^{\infty} \frac{e(\pm\gamma)}{\sqrt{2\gamma}}
d\gamma.$$
A cognate argument shows that also
$$\int_{-\mu}^0
e\big(qjF(\alpha)\big)d\alpha=\frac{W_{{\rm{sgn}}(j)}}{\sqrt{qc_2|j|}}
+ O\big((q|j|)^{(-1-\theta)/2}\big).$$
When $c_2<0$ perhaps the simplest thing is to observe that this case is
formally equivalent to taking complex conjugates.  Thus, in general, we have
$$\int_{-\mu}^{\mu}
e\big(qjF(\alpha)\big)d\alpha=\frac{2W_{{\rm{sgn}}(c_2j)}}{\sqrt{q|c_2j|}}
+ O\big((q|j|)^{(-1-\theta)/2}\big).$$
Hence
\begin{align*}
\sum_{Q<q\le 2Q} q \int_{\beta_h-\mu}^{\beta_h+\mu}&
e\big(q(jf(\beta)-h\beta)\big)d\beta \\ &=  \sum_{Q<q\le 2Q}
q^{\frac12} e(q\phi_h) \frac{2W_{{\rm{sgn}}(c_2j)}}{\sqrt{|c_2j|}}
+ O\big(Q^{(3-\theta)/2}|j|^{(-1-\theta)/2}\big).\\
\end{align*}
Thus
$$\sum_{Q<q\le 2Q} q \int_{\beta_h-\mu}^{\beta_h+\mu}
e\big(q(jf(\beta)-h\beta)\big)d\beta \ll
Q^{\frac12}\lambda_h^{-1}|j|^{-\frac12} +
Q^{(3-\theta)/2}|j|^{(-1-\theta)/2} \ . $$
Hence, by (\ref{e:two6}),
$$N_5\ll J^{\frac12}Q^{\frac12+\varepsilon} +
J^{-\frac12}(\log J)Q^{\frac32} + J^{(1-\theta)/2}
Q^{(3-\theta)/2}.$$
Thus we have established that
$$\widetilde{N}(Q,\delta)
\ll \delta Q^2 + \delta^{-\frac12}Q^{\frac12+\varepsilon} +
\delta^{\frac12}Q^{\frac32}\log\textstyle{\frac1\delta}
+ Q\log\textstyle{\frac1\delta} +
\delta^{\frac{\theta-1}2}Q^{\frac{3-\theta}2}.$$
When $\frac1\delta\le Q^{1-\varepsilon}\log Q$ we have
$$\delta^{\frac12}\big(\log\textstyle{\frac1\delta}\big)Q^{\frac32}
\le \delta Q^{2-\frac12\varepsilon}(\log Q)^{\frac32}\ll \delta
Q^2$$ and when $\frac1\delta> Q^{1-\varepsilon}\log Q$ we have
$$\delta^{\frac12}\big(\log\textstyle{\frac1\delta}\big)Q^{\frac32}
\ll \delta^{-\frac12}Q^{\frac12+\varepsilon}.$$ Moreover, when
$\frac1\delta\le Q^{1-2\varepsilon}\log^2Q$ we have
$$\big(\log\textstyle{\frac1\delta}\big)Q\ll \delta Q^2$$ and when
$\frac1\delta>Q^{1-2\varepsilon}\log^2Q$ we have
$$\big(\log\textstyle{\frac1\delta}\big)Q\ll
\delta^{-\frac12}Q^{\frac12+\varepsilon}.$$ Therefore
$$\widetilde{N}(Q,\delta) \ll \delta Q^2 +
\delta^{-\frac12}Q^{\frac12+\varepsilon} +
\delta^{\frac{\theta-1}2}Q^{\frac{3-\theta}2}.$$

\noindent This completes the proof of Theorem \ref{t:thm4}.

\Section{The proof of Theorem \ref{t:thm5}} \noindent This is
easily deduced from Theorem \ref{t:thm4} in the same manner that
Theorem \ref{t:thm2} is deduced from Theorem \ref{t:thm1}.

\Section{The proof of Theorem \ref{t:thm3}} \label{secthm3}

\noindent We are given that $ \mathcal{C} $ is a $C^{(2)}$
non-degenerate planar curve. Thus,  $ \mathcal{C} = \mathcal{C}_f
:= \{ (x, f(x)) \in {\mathbb R}^2 : x \in I \}$ for some interval
$I$  of $ {\mathbb R} $ and $f \in C^{(2)}(I)$. 
Also,
since $ \mathcal{C}_f$ is non-degenerate we have that  $f''(x)
\neq 0 $ for almost all $x \in I $. Throughout, $\psi$ is an
approximating function such that
$$ {\textstyle \sum_{t=1}^{\infty} } \psi(t)^2  \ < \ \infty  \ . $$

The claim is that $|{\mathcal C}_f \cap {\mathcal
S}(\psi)|_{{\mathcal C}_f}= 0 $.

\noindent{\bf Step 1. \ } We show that there is no loss of
generality in assuming that
\begin{equation}
\psi (t) \ \geq \ t^{-\frac12} (\log t )^{-1} \quad {\rm for \ all
\ } t  . \label{s1}
\end{equation}
To this end, define $\Psi : t \to  \Psi(t)  := \max \{ \psi(t), \,
t^{-\frac12}  (\log t)^{-1} \}  $.
 Clearly, $\Psi$ is an approximating
function and furthermore $ \sum \Psi(t)^2 <\infty $. By
definition,  $\mathcal{S}(\psi) \subset \mathcal{S}(\Psi)$ and so
it suffices to establish the claim with $\psi$ replaced by $\Psi$.
Hence, without loss of generality, (\ref{s1}) can be assumed.

\noindent{\bf Step 2. \ } Let  $\Omega_{f,\psi}$  be the set of
$x\in I$ such that the system of inequalities
\begin{equation}\label{s2}
\left\{
\begin{array}{l}
  \big|x-\frac{p_1}{q}\big|<\frac{\psi(q)}{q}  \\[1ex]
  \big|f(x)-\frac{p_2}{q}\big|<\frac{\psi(q)}{q}
\end{array}
\right.,
\end{equation}
is satisfied for infinitely many ${\vv p}/q\in\Q^2$ with $p_1/q\in
I$. Notice that since $f$ is continuously differentiable, the map
$x\mapsto(x,f(x))$ is locally bi-Lipshitz and so $$ |{\mathcal
C}_f \cap {\mathcal S}(\psi)|_{{\mathcal C}_f}= 0 \quad
\Longleftrightarrow \quad  |\Omega_{f,\psi}|_{\R}= 0 \ .
$$
Hence, it suffices to show that
\begin{equation}\label{s3} |\Omega_{f,\psi}|_{\R} = 0 \ .
\end{equation}

\noindent{\bf Step 3. \ } Next, without loss of generality, we can
assume that $I$ is open in $\R$. Notice that the set $B:=\{ x \in
I: |f''(x)| = 0 \}$ is closed in $I$.  Thus the set  $G:=
I\setminus B $ is open and a standard argument allows one to write
$G$ as a countable union of bounded intervals $I_i$ on which $f$
satisfies
\begin{equation}
\label{infsup} 0\ < \ c_1 \ := \ \inf_{x \in I_0} |f''(x)|  \ \leq
\ c_2 \ := \ \sup_{x \in I_0} |f''(x)| \ < \ \infty  \ .
\end{equation}
The constants $c_1, c_2 $ depend on the particular choice of
interval $I_i$. For the moment, assume that $|\Omega_{f,\psi} \cap
I_i |_{\R} = 0 $ for any $i \in {\mathbb N}$. On using the fact
that $|B|_{\R} = 0 $, we have that
$$ |\Omega_{f,\psi}|_{\R}  \ \leq \  | B \cup \mbox{ \small
$\bigcup\limits_{i=1}^{\infty}$ } (\Omega_{f,\psi}\cap I_i)|_{\R}
\ \leq \ | B |_{\R} \, + \,  \mbox{ \small
$\sum\limits_{i=1}^{\infty}$ } |\Omega_{f,\psi}\cap I_i|_{\R} = 0
$$
and this establishes (\ref{s3}). Thus, without loss of generality,
and  for the sake of clarity we assume that $f$ satisfies
(\ref{infsup}) on $I$ and that $I$ is bounded. The upshot of this
is that $f$ satisfies the conditions imposed in Theorem
\ref{t:thm1}.

\noindent{\bf Step 4. \ } For a point ${\vv p}/q\in\Q^2$, denote
by $\sigma(\vv p/q)$ the set of $x\in I$ satisfying (\ref{s2}).
Trivially, \begin{equation} \label{s3a} |\sigma(\vv p/q)|_{\R} \,
\le \, 2 \psi(q)/q  \ . \end{equation}

\noindent  Assume that $\sigma(\vv p/q)\not=\emptyset$ and let $x
\in \sigma(\vv p/q)$. By the mean value theorem,
$f(x)=f(p_1/q)+f'(\tilde x)(x-p_1/q)$ for some $\tilde x \in I$.
We can assume that $f'$ is bounded on $I$ since $f''$ is bounded
and $I$ is a  bounded interval. Suppose $2^n\le q< 2^{n+1}$. By
(\ref{s2}),
$$\textstyle{ \big|f(\frac{p_1}{q})-\frac{p_2}{q}\big| \ \le \
\big|f(x)-\frac{p_2}{q}\big| \ + \ \big|f'(\tilde
x)\big(x-\frac{p_1}{q}\big)\big| \ < \  c_3 \, \psi(q)/q \ \le \
c_3 \, \psi(2^n)/2^n}  \ , $$ where  $c_3>0$ is a constant. Thus,
\begin{eqnarray*}
\mathrm{card} \{{\vv p}/q\in\Q^2 \!\!\!\! &:& \!\!\!\! 2^n\le q<
2^{n+1}, \, \sigma(\vv p/q)\not=\emptyset \}   \\ & \leq &
\mathrm{card} \left\{{\vv p}/q\in < \Q^2: q \leq 2^{n+1}, \, p_1/q
\in I,  \, \textstyle{ \big|f(\frac{p_1}{q})-\frac{p_2}{q}\big| <
c_3
\, \psi(2^n)/2^n} \right\} \\
& \leq & \mathrm{card} \left\{a/q\in\mathbb Q : q \leq 2^{n+1}, \,
a/q \in I, \,  \textstyle{ \big\|q f(\frac{a}{q})\big\| < 2 c_3 \,
\psi(2^n) } \right\} \ .
\end{eqnarray*}
In view of (\ref{s1}), Theorem \ref{t:thm1} implies that
\begin{equation}
\label{s4} \mathrm{card} \{{\vv p}/q\in\Q^2 : 2^n\le q< 2^{n+1},
\, \sigma(\vv p/q)\not=\emptyset \} \ll \psi(2^n) \, 2^{2n} \ .
\end{equation}

\noindent{\bf Step 5. \ } For $n \geq 0$, let
$$
\Omega_{f,\psi}(n) \ :=  \bigcup_{\vv p/q\in\Q^2,\,\sigma(\vv p/q)
\not=\emptyset,\,2^n\le q<2^{n+1}}
\!\!\!\!\!\!\!\!\!\!\!\!\!\!\!\!\! \sigma(\vv p/q)  \  \  .
$$

\noindent Then  $|\Omega_{f,\psi}|_{\R} = |\limsup_{n \to \infty}
\Omega_{f,\psi}(n)|_{\R} $  and the  Borel-Cantelli Lemma implies
(\ref{s3}) if   $\textstyle{ \sum_{n=0}^{\infty} } |
\Omega_{f,\psi}(n)|_{\R} < \infty $.
In view of (\ref{s3a}) and (\ref{s4}), it follows
that
\begin{eqnarray*}
\sum_{n=0}^{\infty}  | \Omega_{f,\psi}(n)|_{\R}
 & = & \sum_{n=0}^\infty \ \ \ \ \sum_{\vv p/q\in\Q^2,\,\sigma(\vv p/q)
\not=\emptyset,\,2^n\le q<2^{n+1}}
\!\!\!\!\!\!\!\!\!\!\!\!\!\!\!\!\!\! |\sigma(\vv p/q)|_{\R}
\\ \\ & \ll &
\sum_{n=0}^\infty \psi(2^n)/ 2^n  \times \psi(2^n) \, 2^{2n} \
\asymp \ \sum_{t=1}^{\infty} \psi(t)^2 \  < \ \infty \ \ .
\end{eqnarray*}

\noindent This completes the proof of Theorem \ref{t:thm3}.


\Section{The proof of Theorem \ref{t:thm6}}

\noindent In spirit, the proof of Theorem \ref{t:thm6} follows the
same line of argument as the proof of Theorem \ref{t:thm3}.
Throughout, $s \in (1/2,1)$ and  $\psi$ is an approximating
function such that $$ {\textstyle \sum_{t=1}^{\infty} }t^{1-s}
\psi(t)^{s+1}  \  < \   \infty  \ . $$

\noindent{\bf Step 1. \ } Choose $\eta > 0 $ such that $ \eta <
(2s-1)/(s+1) $. Note that $(2s-1)/(s+1) $ is strictly positive
since $s > 1/2$. By considering the auxiliary function $\Psi : t
\to  \Psi(t)  := \max \{ \psi(t), \, t^{-1+\eta} \}  $,  it is
easily verified that there is no loss of generality in assuming
that
\begin{equation}
\psi (t) \ \geq \ t^{-1+\eta}  \quad {\rm for \ all \ } t  .
\label{s11}
\end{equation}

\noindent{\bf Step 2. \ } Let  $\Omega_{f,\psi}$  be defined via
the  system of inequalities (\ref{s2}) as in Step 2 of
\S\ref{secthm3}. On making use of the fact that  the map
$x\mapsto(x,f(x))$ is locally bi-Lipshitz  we have that  $$
\hs({\mathcal C}_f \cap {\mathcal S}(\psi))= 0 \quad
\Longleftrightarrow \quad \hs(\Omega_{f,\psi})= 0 \ . $$ Hence, it
suffices to show that $ \hs(\Omega_{f,\psi}) = 0$.

\noindent{\bf Step 3. \ }Let $B:=\{ x \in I: |f''(x)| = 0 \}$.
Since $\dim B \leq 1/2$ and $s > 1/2$, it follows from the
definition of $\hs$ that $\hs(B) =0$. As in Step 3 of
\S\ref{secthm3}, the set $G:= I\setminus B $ can be  written as a
countable union of bounded intervals $I_i$ on which $f$ satisfies
(\ref{infsup})   and moreover we can assume that  $|I_i |_{\R}
\leq 1$. Since $f \in C^{(3)} (I) $, it follows that  $ | f''(x) -
f''(y) | \ll | x - y | \leq | x - y |^{\theta} $ for any $x,y \in
I_i $ and $ 0 \leq \theta \leq 1$; i.e. $f'' \in
{\rm{Lip}}_{\theta}(I_i)$.  In particular, with Theorem
\ref{t:thm4} in mind, we may take
$$ 1> \theta > \textstyle{\frac{2- 3\eta}{2-\eta}} \ . $$
Now the same argument as in Step 3 of \S\ref{secthm3} with
Lebesgue measure $ |\ \ |_{\R}$ replaced by Hausdorff measure
$\hs$, enables us to conclude that $f$ satisfies $(\ref{infsup}) $
on $I$ and moreover the conditions imposed in Theorem \ref{t:thm4}
are satisfied.

\noindent{\bf Step 4. \ } This is exactly as in Step 4 of
\S\ref{secthm3} apart from the fact  that the conclusion
(\ref{s4}) follows as a consequence of (\ref{s11}) and Theorem
\ref{t:thm4}.

\noindent{\bf Step 5. \ } With $\Omega_{f,\psi}(n)$ as in
Step 5 of \S\ref{secthm3}, we have that for each $l \in  {\mathbb
N}$, $$ \{  \Omega_{f,\psi}(n) : n=l, \, l+1, \ldots  \, \} $$ is
a cover for $ \Omega_{f,\psi}$ by sets $\sigma(\vv p/q)$ of
maximal diameter  $2 \psi(2^l)/2^l $. This  makes use of the
trivial fact that each set $\sigma(\vv p/q)$ is contained in an
interval of length at most $2 \psi(q)/q $.  It follows from the
definition of  Hausdorff measure that with $\rho := 2
\psi(2^l)/2^l $,
\begin{eqnarray*}
\hs_{\rho} ( \Omega_{f,\psi} ) & \leq  & \sum_{n=l}^{\infty}  \ \
\ \sum_{\vv p/q\in\Q^2,\,\sigma(\vv p/q) \not=\emptyset,\,2^n\le
q<2^{n+1}} \!\!\!\!\!\!\!\!\!\!\!\!\!\!\!\!\!\! ( 2\, \psi(2^n)/
2^n)^s
\\ \\ & \ll & \sum_{n=l}^{\infty} \
 (\psi(2^n)/ 2^n)^s   \times \psi(2^n) \, 2^{2n}  \ \longrightarrow \ 0 \ \
\end{eqnarray*}
 as $\rho \to 0$; or equivalently at $l \to \infty$.   Hence,
 $ \hs ( \Omega_{f,\psi} ) = 0 $ and this completes the
 proof of Theorem \ref{t:thm6}.

\Section{Various generalizations: the multiplicative setup}

For the sake of brevity, we shall restrict our  attention to  the
Lebesgue theory only.

\noindent Given approximating functions $\psi_1,\psi_2$, a point
$\vv y\in\R^2$ is said to be  {\it simultaneously
$(\psi_1,\psi_2)$--approximable} if there are infinitely many
$q\in {\mathbb N}$ such that $$ \|q y_i\|<\psi_i(q) \hspace{17mm}
1\le i\le 2 \ .  $$ Let  ${\mathcal S}(\psi_1,\psi_2)$  denote the
set of simultaneously $(\psi_1,\psi_2)$--approximable points in
$\R^2$. This set 
is clearly a generalization of  ${\mathcal S}(\psi)$ in which
$\psi = \psi_1 = \psi_2$. The following statement is a natural
generalization of Khinchin's theorem: $$|{\mathcal
S}(\psi_1,\psi_2)|_{\R^2} =\left\{\begin{array}{ll} \mbox{\rm
Z{\scriptsize ERO}} & {\rm if} \;\;\; \sum \;  \psi_1(t)\,
\psi_2(t) \;\; <\infty  \\ & \\ \mbox{\rm F{\scriptsize ULL}} &
{\rm if} \;\;\; \sum  \;  \psi_1(t)\, \psi_2(t) \;\;
 =\infty \; \;
\end{array}\right..$$

Next, given an approximating function $\psi$, a point $\vv y
\in\R^2$ is said to be {\it multiplicatively $\psi$--approximable}
if there are infinitely many $q\in {\mathbb N}$ such that $$
\textstyle{\prod_{i=1}^2} \|qy_i\| \ < \ \psi(q) \ .
  $$ Let ${\mathcal S}^*(\psi)$ denote the set of multiplicatively
$\psi$--approximable points in $\R^2$. In view of  Gallagher's
theorem we have that: $$|{\mathcal S}^*(\psi)|_{\R^2}
=\left\{\begin{array}{ll} \mbox{\rm Z{\scriptsize ERO}} & {\rm if}
\;\;\; \sum \;  \psi(t)^2 \ \log t \;\; < \ \infty\\ & \\
\mbox{\rm F{\scriptsize ULL}} & {\rm if} \;\;\; \sum  \;
\psi(t)^2 \ \log t  \;\;
 = \ \infty \; \;
\end{array}\right..$$

Now  let $\mathcal{C}$ be a $C^{(3)}$ non-degenerate planar curve.
The goal is to obtain the analogues of the above `zero-full'
statements for the sets  $\mathcal{C} \cap {\mathcal
S}(\psi_1,\psi_2)$ and $ \mathcal{C} \cap {\mathcal S}^*(\psi)$.
It is highly likely that the counting results obtained in this
paper, in particular Theorem \ref{t:thm4}, together with the ideas
developed in \cite{BV} will yield the following convergence
statements.

\noindent {\bf Claim 1. \ } $ |\mathcal{C} \cap {\mathcal
S}(\psi_1,\psi_2) |_{\mathcal{C}} = 0  \ \ $ if $ \ \
\sum\psi_1(t) \psi_2(t) < \infty $.

\noindent {\bf Claim 2. \ } $ |\mathcal{C} \cap {\mathcal
S}^*(\psi) |_{\mathcal{C}} = 0  \ \ $ if $ \ \ \sum\psi(t) \log t
< \infty $.

In the case that the planar curve $\mathcal{C}$ belongs to a
special class of  rational quadrics, both these claims have been
established in \cite{BV}. Furthermore, in \cite{BV} the divergent
analogue of Claim 1 has been established. Thus, establishing Claim
1 would complete the Lebesgue theory for simultaneously
$(\psi_1,\psi_2)$--approximable points on planar curves.

Currently, D. Badziahin is attempting to establish the above
claims and is also investigating the Hausdorff measure theory.

\vspace{3ex}

\noindent{\em Acknowledgements: } SV would like to thank Victor
Beresnevich for the numerous enlightening conversations regarding
the general area of Diophantine approximation on manifolds and for
so generously sharing his insight. He would also like to thank
those simply wonderful  girls Ayesha and Iona for introducing him
to the `fourth dimension'  at our very  special place -- Almscliff
Crag.


\end{document}